\documentclass[11pt]{article}
\usepackage[dvips]{graphics}
\input{amssym.def}
\input{amssym.tex}

\title{\Large\bf  Helmholtz equation in a semi-infinite strip with impedance boundary conditions
of the third and fifth orders}
\author{\bf {\sc By Y.A.\ Antipov}\\ 
Department of Mathematics, Louisiana State University\\
Baton Rouge LA 70803, USA}

\date{}

\newcommand{\I}{\mathop{\rm Im}\nolimits}
\newcommand{\R}{
\mathop{\rm Re}\nolimits}
\newcommand{\const}{\mbox{const}}

\newcommand{\Md}{\partial}

\newcommand{\Ga}{\alpha}
\newcommand{\Gb}{\beta}
\newcommand{\Gd}{\delta}

\newcommand{\Gve}{\varepsilon}
\newcommand{\Gf}{\phi}

\newcommand{\Gg}{\gamma}

\newcommand{\Gk}{\kappa}

\newcommand{\Gl}{\lambda}
\newcommand{\Gn}{\eta}
\newcommand{\Gm}{\mu}

\newcommand{\Gt}{\theta}

\newcommand{\Gr}{\rho}

\newcommand{\Gs}{\sigma}

\newcommand{\Go}{\omega}
\newcommand{\Gx}{\xi}
\newcommand{\Gy}{\psi}
\newcommand{\Gz}{\zeta}
\newcommand{\GD}{\Delta}
\newcommand{\GF}{\Phi}

\newcommand{\GL}{\Lambda}

\newcommand{\CJ}{{\cal J}}

\newcommand{\beq}{\begin{equation}}
\newcommand{\eeq}{\end{equation}}
\newcommand{\barr}{\begin{eqnarray}}
\newcommand{\earr}{\end{eqnarray}}
\newcommand{\beqn}{\begin{equation*}}
\newcommand{\eeqn}{\end{equation*}}
\newcommand{\barrn}{\begin{eqnarray*}}
\newcommand{\earrn}{\end{eqnarray*}}

\newcommand{\fr}{\frac}

\begin{document}
\maketitle

\begin{abstract}

Two boundary value problems for the Helmholtz equation in a semi-infinite strip  are considered.
The main feature of these problems is that, in addition to the function and its normal derivative
on the boundary, the functionals of the boundary conditions possess tangential derivatives of
the second and fourth orders. Also, the setting of the problems is complimented by certain edge
conditions at the two vertices of the semi-strip.  The problems model wave propagation in a 
semi-infinite waveguide with membrane and plate walls.
A technique for the exact solution of these fluid-structure interaction problems is proposed. It requires application of 
two Laplace transforms with respect to both variables with the parameter of the second transform being
a certain function  of the first Laplace transform parameter. Ultimately, this method yields two scalar Riemann-Hilbert
problems with the same coefficient and different right-hand sides. The dependence of the existence and uniqueness results of the physical model problems 
 on the index of the Riemann-Hilbert problem is discussed.

\end{abstract}

\setcounter{equation}{0}

\section{Introduction}

Boundary value problems for the Helmholtz equation with order $n\ge 2$ derivatives in the boundary conditions have been employed in the theory of diffraction since the work \cite{ryt}, where the second order derivatives
on the boundary were used to model the surfaces of highly conducted materials. The first
order impedance boundary conditions were generalized in \cite{kar} by adding second order tangential 
derivatives on the surface in order to model metal-backed dielectric layers. Order
$n\ge 3$  boundary and transition conditions in electromagnetic diffraction theory were systematically 
studied in \cite{sen}.

Higher order tangential derivatives in the boundary conditions naturally arise in model problems of aerodynamic noise theory and underwater acoustics when sound waves in fluids interact with flexible surfaces  of waveguides
\cite{jun}.    
Exact solutions and their analysis are available \cite{lep} for problems on a compressible fluid bounded by an
infinite membrane and  an elastic plate fixed along two or more parallel lines when the system is excited by an incident plane wave. These models for membranes and elastic plates are governed by the Helmholtz equation with
 the third and fifth order derivatives in the boundary conditions, respectively. 
The Wiener-Hopf method was applied   in  \cite{kou1} to study the motion of an infinite plane composed of 
two half-planes with different 
elastic constants due to hydroacoustic pressure in the fluid beneath the plate. 
Edge diffraction of an incident acoustic plane wave by a thin elastic half-plane was treated in  \cite{can}. 
Due to the geometry of the model problem it was also solved by the Wiener-Hopf method. 

For more complicated domains like  a wedge, two joint wedges, or a semi-infinite strip whose
boundaries are composed of either  membranes or elastic plates the Wiener-Hopf method
is not applicable. The Buchwald method \cite{buc} proposed for the solution of the model problem on
diffraction of Kelvin waves at a corner was further developed \cite{kou2} to study diffraction of acoustic waves 
in a semi-infinite waveguide whose surface is formed by elastic plates.
In this work the authors determined that the number of free constants
to be determined from the conditions at the vertices of the structure depends only on the orders of the 
derivatives in the boundary conditions. They also found  an integral representation of the acoustic pressure distribution.
Some model problems of sound-structure interaction with high-order boundary conditions by the method of eigenfunction expansions 
were treated in \cite{law}. These include the problem of propagation of an acoustic wave in a semi-infinite
waveguide when the upper boundary is a semi-infinite membrane,  while the lower and the finite vertical sides
are acoustically rigid walls (in this case, by symmetry,  the problem for a a half-strip
reduces to a problem for an infinite strip).

The Poincar\'e boundary value problem for the modified Helmholtz operator $\GD-k^2$ ($k$ is a real number)
in a semi-infinite strip
was studied in \cite{ant1} by the method proposed in \cite{fok}. It was shown that the problem
reduces to an order-2 vector Riemann-Hilbert problem on the real axis whose matrix coefficient, in general,
does not admit an explicit factorization by the methods currently available in the literature. However,
in some important cases, including the case of impedance boundary conditions,
it allows for an exact factorization and in some other cases the vector problem reduces to a problem with triangle matrix coefficient, or could be even decoupled. Notice that  in the case of
the Helmholtz operator $\GD+k^2$ ($k\in{\Bbb R}$), the contour of the Riemann-Hilbert problem comprises two semi-infinite
rays and two circular arcs.

Our goal in this work is to develop an efficient method
for the Helmholtz equation in a  semi-infinite strip $\{0<x<\infty, 0<y<a\}$
with generalized impedance boundary conditions
of higher order. The main feature of this technique is that it applies  two Laplace transforms in
a nonstandard  way. The first  transform is utilized with respect to $x$ in the classical way, while the second one 
is applied in the $y$-direction (the function is extended by zero for $y>a$)  with a parameter $\Gz$ that is
a root of the characteristic polynomial of the ordinary
differential operator $d^2/dy^2+k^2-\Gn^2$, the Laplace image of the original Helmholtz operator. 
This root is $\Gz=\sqrt{\Gn^2-k^2}$, where $k$
is the wave number, $\Gn$ is the parameter of the first transform, and $\sqrt{\Gn^2-k^2}$ is a fixed
branch of the function $\Gz^2=\Gn^2-k^2$.  The method to be proposed ultimately yields
two symmetric scalar Riemann-Hilbert problems on the real axis equivalent to the original model problem.
We emphasize that the contour is the real axis regardless if $k$ is real, imaginary, or a general complex number. 
The Riemann-Hilbert problems share the same coefficient and have different right-hand sides.
Remarkably, the coefficient is a simple rational function, $G(\Gn)=P_n(\Gn)/P_n(-\Gn)$, where the degree
of the polynomial $P_n(\Gn)$ is $n=3$ (in the membrane case) and $n=5$ (in the elastic plate case)
and coincide with the order of the highest derivative
involved in the boundary conditions. We determine the number of free constants in the solution that depends not only
on the order of the tangential derivatives in the boundary conditions but also on the position
of the zeros of the polynomial $P_n(\Gn)$ and therefore the parameters of the problem. Also, 
we derive explicitly a system of linear algebraic equations for the unknown constants and summarize the results
by stating an existence - uniqueness theorem. Finally, we write down representation formulas for the solution by quadratures and, in addition,  by series  convenient for computational purposes.

\setcounter{equation}{0}

\section{Helmholtz equation in a semi-infinite waveguide: membrane walls}\label{helm}

\subsection{Formulation}\label{form}

Of concern is the Helmholtz equation
\beq
\left(\fr{\Md^2}{\Md x^2}+\fr{\Md^2}{\Md y^2}+k^2\right)u(x,y)=g(x,y), \quad 0<x<\infty, \quad 0<y<a,
\label{2.1}
\eeq
with respect to an unknown function $u(x,y)$ subjected to the boundary conditions 
$$
\left[\left(\fr{\Md^2}{\Md x^2}+\Ga_0^2\right)\fr{\Md}{\Md y}-\mu_0\right]u=g_0(x), \quad  (x,y)\in W_0=\{0<x<\infty, \; y=0\},
$$
$$
\left[\left(\fr{\Md^2}{\Md x^2}+\Ga_1^2\right)\fr{\Md}{\Md y}+\mu_1\right]u=g_1(x), \quad  (x,y)\in W_1=\{0<x<\infty, \quad y=a\},
$$
\beq
\left[\left(\fr{\Md^2}{\Md y^2}+\Ga_2^2\right)\fr{\Md}{\Md x}-\mu_2\right]u=g_2(y), \quad  (x,y)\in W_2=\{ x=0, \quad 0<y<a\}.
\label{2.2}
\eeq
This boundary value problem governs acoustic wave propagation   in a semi-infinite waveguide (Fig.1).
\begin{figure}[t]
\centerline{
\scalebox{0.6}{\includegraphics{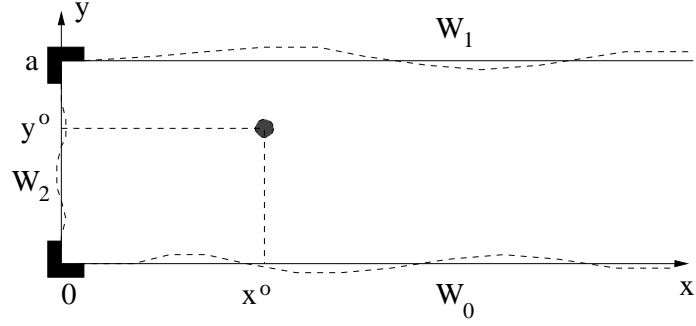}}}
\caption{Geometry of the problem on a semi-infinite waveguide.}
\label{fig1}
\end{figure} 
Here,
$\R[e^{-i\Go t}u(x,y)]$ is the fluid velocity potential, $\Go$ is the frequency, $\Go=\Go_1+i\Go_2$,
$\Go_j>0$, $t$ is time, $k=\Go/c$ is
the wave number, $c$ is the sound speed in the fluid. The pressure distribution $p(x,y)$ and the
deflection of the horizontal and vertical walls, $u_0(x)$, $u_1(x)$, and $u_2(y)$, are expressed through the velocity potential as
\beq
p(x,y)=i\Go\Gr u(x,y), \quad u_j(x)=\fr{i}{\Go}u_y(x,y_j), \quad u_2(y)=\fr{i}{\Go}u_x(0,y),  
\label{2.3}
\eeq 
where $j=0,1$,  $y_0=0$,  $y_1=a$,  
$\Gr$ is the mean fluid density, and 
 the suffixes $x$  and $y$ denote differentiation with respect to $x$ and $y$,
respectively,
The boundary conditions (\ref{2.2}) model the deflection of the membrane walls
due to pressure loading (Leppington, 1978). The parameters involved are
\beq
\Ga_j=\Go\sqrt{\fr{m_j}{T_j}}, \quad \mu_j=\fr{\Gr\Go^2}{T_j},\quad j=0,1,2,
\label{2.4}
\eeq
where $m_j$  is the mass per unit area, and  $T_j$ is the surface tension for the membrane $W_j$. Since 
$\I\Go>0$, we have $\I\Ga_j>0$.
The functions $g(x,y)$, $g_0(x)$, $g_1(x)$, and $g_2(y)$ are prescribed. They will be later selected as $g_0(x)=g_1(x)=0$,
$0<x<\infty$, $g_2(y)=0$, $0<y<a$, and $g(x,y)=-\Gd(x-x^\circ)\Gd(y-y^\circ)$, $(x^\circ,y^\circ)$ is an internal point of the semi-strip,
and $\Gd(\cdot)$ is the Dirac function.
We also need to specify the edge conditions at $x=y=0$ and $x=0$, $y=a$. It is assumed that the edges are fixed, and therefore
the following four conditions have to be satisfied:
$$
\lim_{x\to 0^+}u_y(x,0)=\lim_{x\to 0^+}u_y(x,a)=0,
$$
\beq
\lim_{y\to 0^+}u_x(0,y)=\lim_{y\to a^-}u_x(0,y)=0.
\label{2.5}
\eeq

\subsection{Two scalar Riemann-Hilbert problems}\label{RHP}

Our goal in this section is to present a method that is capable to convert the boundary value problem for the Helmholtz
equation (\ref{2.1})
with higher-order boundary conditions (\ref{2.2}) in a semi-strip to two scalar Riemann-Hilbert problems. To achieve this, we apply first
the Laplace transform with respect to $x$
\beq
\tilde u(\Gn,y)=\int_0^\infty u(x,y)e^{i\Gn x}dx
\label{2.5'}
\eeq
to equation (\ref{2.1}) and the first and second boundary conditions in (\ref{2.2}). In conjunction with (\ref{2.5})
this brings us to the one-dimensional boundary value problem
$$
L[\tilde u]\equiv \left(\fr{d^2}{dy^2}-\Gz^2\right)\tilde u(\Gn,y)=f(y), \quad 0<y<a,
$$$$
U_0[\tilde u]\equiv -\tilde u_y(\Gn,0)+\tilde \mu_0(\Gn)\tilde u(\Gn,0)=\tilde g^0(\Gn),
$$
\beq
U_1[\tilde u]\equiv\tilde u_y(\Gn,a)+\tilde \mu_1(\Gn)\tilde u(\Gn,a)=\tilde g^1(\Gn),
\label{2.6}
\eeq
where 
$$
\Gz^2=\Gn^2-k^2, \quad 
f(y)=u_x(0,y)-i\Gn u(0,y)+\tilde g(\Gn,y),
$$
$$
\tilde\mu_j=\fr{\mu_j}{\Ga_j^2-\Gn^2}, \quad \tilde g^j(\Gn)=(-1)^{j+1}\fr{\tilde g_j(\Gn)+c_j}{\Ga_j^2-\Gn^2},
$$
\beq
\tilde g(\Gn,y)=\int_0^\infty g(x,y)e^{i\Gn x}dx, \quad \tilde g_j(\Gn)=\int_0^\infty g_j(x)e^{i\Gn x}dx, \quad j=0,1,
\label{2.8}
\eeq
and $c_0$ and $c_1$ are unknown constants 
\beq
c_0=u_{xy}(0^+,0), \quad c_1=u_{xy}(0^+,a).
\label{2.9}
\eeq
These constants and the functions $u_x(0,y)$ and $u(0,y)$ in (\ref{2.8}) are to be determined.

Denote by  $\Gz=\sqrt{\Gn^2-k^2}$   the single branch of the two-valued function $\Gz^2=\Gn^2-k^2$
in the $\Gn$-plane cut along the straight line joining the branch points $\Gn=\pm k$ ($k\in{\Bbb C}^+$) and passing
through the infinite point; the branch is fixed by the condition $\Gz=-ik$ as $\Gn=0$. We next  employ the Green function
of the boundary value problem (\ref{2.6})
$$
G(y,s)=-\fr{e^{-\Gz|y-s|}}{2\Gz}+\fr{1}{2\Gz \tilde\GD(\Gz)}\left\{(\tilde\mu_0-\Gz)
[\Gz\cosh\Gz(a-y)+\tilde\mu_1\sinh\Gz(a-y)]e^{-\Gz s}\right.
$$
\beq
\left.
+(\tilde\mu_1-\Gz)(\Gz\cosh \Gz y+\tilde\mu_0\sinh \Gz y)e^{-\Gz(a-s)}\right\},
\label{2.10}
\eeq
and the fundamental system, namely the two solutions $\Gf_0(y)$ and $\Gf_1(y)$ of the problem
$$
L[\Gf_j(y)]=0, \quad 0<y<a,
$$
\beq
U_m[\Gf_j]=\Gd_{mj}, \quad m,j=0,1,
\label{2.11}
\eeq
which are
\beq
\Gf_0(y)=\fr{\Gz\cosh\Gz(a-y)+\tilde\mu_1\sinh\Gz(a-y)}{\tilde \GD(\Gz)},\quad
\Gf_1(y)=\fr{\Gz\cosh\Gz y+\tilde\mu_0\sinh\Gz y}{\tilde \GD(\Gz)}.
\label{2.12}
\eeq
Here,
\beq
\tilde\GD(\Gz)=(\tilde\mu_0+\tilde\mu_1)\Gz\cosh a\Gz +(\tilde\mu_0\tilde\mu_1+\Gz^2)\sinh a\Gz .
\label{2.13}
\eeq
In terms of the Green function and the functions $\Gf_0$ and $\Gf_1$ the solution of the problem
(\ref{2.6}) can be written as
\beq
\tilde u(\Gn,y)=\int_0^a G(y,s)f(s)ds+\tilde g^0(\Gn)\Gf_0(y)+\tilde g^1(\Gn)\Gf_1(y).
\label{2.14}
\eeq
Now, by putting $y=0$ and $y=a$ in this formula and denoting 
$$
\hat u(0,i\Gz)=\int_0^a u(0,y)e^{-\Gz y}dy, \quad \hat u_x(0,i\Gz)=\int_0^a u_x(0,y)e^{-\Gz y}dy,
$$
\beq
\hat{\tilde g}(\Gn, i\Gz)=\int_0^a\tilde g(\Gn,y)e^{-\Gz y}dy,
\label{2.15}
\eeq
we obtain two relations which can be complimented by
their counterparts  with $\Gn$ being replaced by $-\Gn$; the four resulting equations are
$$
\tilde u(\pm \Gn,0)+\GL_{00}(\Gz)\hat u_x(0,i\Gz)\mp i\Gn\GL_{00}(\Gz)\hat u(0,i\Gz)
$$$$
+\GL_{01}(\Gz)\hat u_x(0,-i\Gz)\mp i\Gn\GL_{01}(\Gz)\hat u(0,-i\Gz)+h_0(\pm\Gn)=0,
$$
$$
\tilde u(\pm \Gn,a)+\GL_{10}(\Gz)\hat u_x(0,i\Gz)\mp i\Gn\GL_{10}(\Gz)\hat u(0,i\Gz)
$$
\beq
+\GL_{11}(\Gz)\hat u_x(0,-i\Gz)\mp i\Gn\GL_{11}(\Gz)\hat u(0,-i\Gz)+h_1(\pm\Gn)=0,
\label{2.16}
\eeq
where
$$
\GL_{00}(\Gz)=\fr{1}{2\Gz}\left[1-\fr{\tilde\mu_0-\Gz}{\tilde\GD(\Gz)}(\Gz\cosh a\Gz +\tilde\mu_1\sinh a\Gz)\right],
$$
$$
\GL_{01}(\Gz)=-\fr{\tilde\mu_1-\Gz}{2\tilde\GD(\Gz)}e^{-a\Gz}, \quad \GL_{10}(\Gz)=-\fr{\tilde\mu_0-\Gz}{2\tilde\GD(\Gz)}, 
$$
\beq
\GL_{11}(\Gz)=\fr{e^{-a\Gz}}{2\Gz}\left[1-\fr{\tilde\mu_1-\Gz}{\tilde\GD(\Gz)}(\Gz\cosh a\Gz +\tilde\mu_0\sinh a\Gz)\right].
\label{2.17}
\eeq
and
$$
h_0(\Gn)=\GL_{00}(\Gz)\hat{\tilde g}(\Gn,i\Gz)+\GL_{01}(\Gz)\hat{\tilde g}(\Gn,-i\Gz)-
\fr{(\Gz\cosh a\Gz+\tilde\mu_1\sinh a\Gz)\tilde g^0(\Gn)}{\tilde\GD(\Gz)}-\fr{\Gz\tilde g^1(\Gn)}{\tilde\GD(\Gz)},
$$
\beq
h_1(\Gn)=\GL_{10}(\Gz)\hat{\tilde g}(\Gn,i\Gz)+\GL_{11}(\Gz)\hat{\tilde g}(\Gn,-i\Gz)-\fr{\Gz \tilde g^0(\Gn)}{\tilde\GD(\Gz)}
-\fr{(\Gz\cosh a\Gz+\tilde\mu_0\sinh a\Gz)\tilde g^1(\Gn).}{\tilde\GD(\Gz)}
\label{2.18}
\eeq
We emphasize two points here. First,  the relative simplicity of the relations (\ref{2.16}) is due to the fact that the functions $\Gz(\Gn)$, $\tilde\mu_0(\Gn)$
and $\tilde\mu_1(\Gn)$ are even, which in turn is a consequence of the absence of odd order derivatives in the 
Helmholtz operator and odd order tangential derivatives in the boundary conditions (\ref{2.2}).
The second point is that the relations (\ref{2.16}) connect the $\pm\Gn$-Laplace transforms of the functions $u(x,0)$  and $u(x,a)$
with the functions $\hat u(0,\pm i\Gz)$ and $\hat u_x(0,\pm i\Gz)$.
If the function $u(0,y)$ and the derivative $u_x(0,y)$ are extended by zero to $y\in(a,\infty)$, then $\hat u(0,\pm i\Gz)$ and
 $\hat u_x(0,\pm i\Gz)$
are their Laplace transforms whose parameters are not arbitrary but functions of $\Gn$, the parameter of the 
first Laplace transform (\ref{2.5'}). These functions  of $\Gn$ are the two zeros $\pm\Gz$ of the characteristic polynomial $r^2-\Gz^2$
of the differential operator $L$ in  (\ref{2.6}).  

We assert that the boundary condition on the vertical side $M_2=\{x=0, 0<y<a\}$ has not been satisfied yet. On extending
the function $g_2(y)$ by zero to the interval $y>a$, applying the Laplace
transform to the third condition in (\ref{2.2}) with the parameters $\Gz$ and $-\Gz$ and utilizing the edge conditions (\ref{2.5})
we discover
\beq
-\hat u_x(0,\pm i\Gz)+\hat\mu_2(\Gz)\hat u(0,\pm i\Gz)=\hat g^2(\pm i\Gz),
\label{2.19}
\eeq
where
$$
\hat\mu_2(\Gz)=\fr{\mu_2}{\Gz^2+\Ga_2^2},\quad \hat g^2(i\Gz)=\fr{-\hat g_2(i\Gz)-c_2+c_3 e^{-a\Gz}}{\Gz^2+\Ga_2^2},
$$
\beq
\hat g_2(i\Gz)=\int_0^a g_2(y)e^{-\Gz y}dy,
\label{2.20}
\eeq
and $c_2$ and $c_3$ are unknown constants, 
\beq
c_2=u_{xy}(0,0^+), \quad c_3=u_{xy}(0,a^-),
\label{2.20.0}
\eeq
to be fixed.

Our intension next is to express the four functions $\hat u(0,\pm i\Gz)$ and $\hat u_x(0,\pm i\Gz)$ through the 
functions $\tilde u(\pm\Gn,0)$ and $\tilde u(\pm\Gn,a)$
from the system (\ref{2.16}) and insert them into the two equations (\ref{2.19}). We have first
\beq
\left(\begin{array}{c}
\hat u_x(0,i\Gz)\\
\hat u(0,i\Gz)\\
\hat u_x(0,-i\Gz)\\
\hat u(0,-i\Gz)\\
\end{array}\right)=\GL(\Gz)
\left(\begin{array}{c}
\tilde u(\Gn,0)+h_1(\Gn)\\
\tilde u(-\Gn,0)+h_1(-\Gn)\\
\tilde u(\Gn,a)+h_2(\Gn)\\
\tilde u(-\Gn,a)+h_2(-\Gn)\\
\end{array}\right),
\label{2.21}
\eeq
where
\beq
\GL(\Gz)=-\fr{1}{2\Gn}\left(\begin{array}{cccc}
\Gn(\tilde\mu_0+\Gz) & \Gn(\tilde\mu_0+\Gz) & \Gn(\tilde\mu_1-\Gz)e^{-a\Gz}  & \Gn(\tilde\mu_1-\Gz)e^{-a\Gz} \\
i(\tilde\mu_0+\Gz) & -i(\tilde\mu_0+\Gz)  & i(\tilde\mu_1-\Gz)e^{-a\Gz}  & -i(\tilde\mu_1-\Gz)e^{-a\Gz}\\
\Gn(\tilde\mu_0-\Gz) & \Gn(\tilde\mu_0-\Gz) & \Gn(\tilde\mu_1+\Gz)e^{a\Gz} & \Gn(\tilde\mu_1+\Gz)e^{a\Gz}\\
i(\tilde\mu_0-\Gz)  & -i(\tilde\mu_0-\Gz)  & i(\tilde\mu_1+\Gz)e^{a\Gz} & -i(\tilde\mu_1+\Gz)e^{a\Gz}\\
\end{array}\right),
\label{2.22}
\eeq
and then, after utilizing equations (\ref{2.19}), we eventually arrive at the following remarkably simple scalar Riemann-Hilbert problems
which share the same coefficient and have different free terms:
\beq
\GF_j^+(\Gn)=H(\Gn)\GF_j(\Gn)+f_j(\Gn), \quad -\infty<\Gn<+\infty,
\label{2.23}
\eeq
subject to the symmetry conditions
\beq
\GF_j^+(\Gn)=\GF_j^-(-\Gn), \quad \Gn\in {\Bbb C}, \quad j=0,1,
\label{2.23'}
\eeq
where 
\beq
\GF_0^\pm(\Gn)=\tilde u(\pm\Gn, 0), \quad \GF_1^\pm(\Gn)=\tilde u(\pm\Gn, a), 
\label{2.24}
\eeq
and
$$
H(\Gn)=-\fr{\Gn+i\hat\mu_2(\Gz)}{\Gn-i\hat\mu_2(\Gz)}=-\fr{\Gn(\Gn^2-k^2+\Ga_2^2)+i\mu_2}{\Gn(\Gn^2-k^2+\Ga_2^2)-i\mu_2},
$$$$
f_j(\Gn)=-h_j(\Gn)+H(\Gn)h_j(-\Gn)+\fr{\Gn}{(\Gn-i\hat\mu_2)\tilde\GD(\Gz)}
$$
\beq
\times\{
[\Gz+(-1)^j\tilde\mu_{1-j}]e^{(1-j)a\Gz}\hat g_2(i\Gz)+[\Gz-(-1)^j\tilde\mu_{1-j}]e^{(j-1)a\Gz}\hat g_2(-i\Gz)\}, \quad j=0,1.
\label{2.25}
\eeq
The functions $\GF_j^+(\Gn)$ and $\GF_j^-(\Gn)$ are analytically continued from the contour, the real axis,
into the upper and lower half-planes, ${\Bbb C}^+$
and ${\Bbb C}^-$, respectively.
Denote 
\beq
q(\Gn)=\Gn(\Gn^2-k^2+\Ga_2^2)+i\mu_2.
\label{2.25.0}
\eeq
Then $H(\Gn)=q(\Gn)/q(-\Gn)$. It turns out
that for all $k=k_1+ik_2$, $k_j\ge 0$, $k_1^2+k_2^2>0$, and all admissible values of the parameters  $\Ga_2$ and $\mu_2$ introduced in
(\ref{2.4}), the zeros  $z_j$  of the cubic polynomial $q(\Gn)$ are simple, and we have the following
three possibilities:

(i) $z_0=-\Gn_0$, $z_1=\Gn_1$, and $z_2=-\Gn_2$, $\I\Gn_j>0$, $j=0,1,2$,

(ii) $z_0=\Gn_0$, $z_1=\Gn_1$, and $z_2=-\Gn_2$, $\I\Gn_0=0$, $\I\Gn_{j}>0$,  $j=1,2$,

(iii) $z_0=\Gn_0$, $z_1=\Gn_1$, and $z_2=-\Gn_2$, $\I\Gn_j>0$, $j=0,1,2$.

\vspace{.1in}

Table 1. The roots $z_j$ $(j=0,1,2)$ of $q(\Gn)=\Gn(\Gn^2-k^2+k^2\Gg_0)+ik^2\Gg_1$ for  some values of the parameters $\Gg_0$, $\Gg_1$, and $k$.

\vspace{1mm}
\begin{tabular}{|c|c|c|c|c|c|}
\hline
$\Gg_0$ &$\Gg_1$ & $k$ & $z_0$ & $z_1$ &$z_2$  \\
\hline
5 &1 & $1+0.1i$ & $-0.0008424 - 0.2540 i$& $-0.2009 + 2.115 i$ & $0.2017 - 1.861 i$\\
\hline
0.5 & 0.1 & $1+0.1i$ &$0.7264 - 0.02424 i$  &$-0.002135 + 0.1872 i$ & $-0.7243 - 0.1629 i$\\
\hline
1 & 0.1 & $1+i$ & $0.5848$ &    $-0.2924 + 0.5065 i$                 &  $-0.2924 - 0.5065 i$ \\
\hline
0.5 &0.05 & $1+0.1i$ &$ 0.7123 + 0.02117 i$ & $-0.0003512 + 0.09816 i$              &$-0.7120 - 0.1193 i$ \\
\hline
1 & 0.1 & $i$ & $-0.4020 + 0.2321 i$ &    $0.4020 + 0.2321 i$                 &  $ - 0.4642 i$ \\
\hline 
\end{tabular}

\vspace{.1in}

In Table 1 we present the roots of the polynomial $q(\Gn)$ in cases (i) (the first two rows),  (ii) (the third row), and
(iii) (the fourth and fifth 
rows).
For convenience, we used the following notations:
$\Gg_0=\Ga_2^2/k^2=m_2c^2/T_2$ and $\Gg_1=\Gm_2/k^2=\Gr c^2/T_2$.

Consider the first case when the polynomial $q(\Ga)$ has one zero in the
upper half-plane and two zeros in the lower half-plane.  
The location of the zeros enables us to factorize the coefficient $H(\Gn)$ as
\beq
H(\Gn)=\fr{H^+(\Gn)}{H^-(\Gn)}, \quad -\infty<\Gn<+\infty, 
\label{2.26}
\eeq
where
\beq
H^+(\Gn)=\fr{(\Gn+\Gn_0)(\Gn+\Gn_2)}{\Gn+\Gn_1}, \quad H^-(\Gn)=-\fr{(\Gn-\Gn_0)(\Gn-\Gn_2)}{\Gn-\Gn_1}, \label{2.27}
\eeq
and $\pm\Gn_j\in {\Bbb C}^\pm$, $ j=0,1,2$.
The index (the winding number) of the function $H(\Gn)$ equals -1, and in the class of functions having
a simple zero at the infinite point one expects  the Riemann-Hilbert problems being solvable if and only if 
a certain condition is fulfilled (Gakhov, 1966). However,  we assert that in our case this condition
is identically satisfied, and the solution always exists. This solution  is unique provided the functions $f_j(\tau)$ are uniquely defined.
Indeed, introduce the Cauchy integrals
\beq
\Psi_j^\pm(\Gn)=\fr{1}{2\pi i}\int_{-\infty}^\infty \fr{f_j(\tau)d\tau}{H^+(\tau)(\tau-\Gn)}, \quad \Gn\in {\Bbb C}^\pm, \quad j=0,1,
\label{2.28}
\eeq
with the densities $f_j(\tau)/H^+(\tau)$ vanishing at the infinite point, and $|f_j(\tau)/H^+(\tau)|\le c|\tau|^{-4}$,
 $\tau\to\pm\infty$, $c$ is a nonzero constant.
The standard application of the continuity principle and the Liouville theorem brings us to the following
representation formulas for the solution of the Riemann-Hilbert problems (\ref{2.25}):
\beq
\GF_j^\pm(\Gn)=H^\pm(\Gn)\Psi_j^\pm(\Gn), \quad 
\Gn\in {\Bbb C}^\pm, \quad j=0,1.
\label{2.29}
\eeq
It is directly verified that
\beq
\fr{f_j(\tau)}{H^+(\tau)}=-\fr{f_j(-\tau)}{H^+(-\tau)},
\label{2.30}
\eeq
and therefore (\ref{2.28}) implies
\beq
\Psi_j^\pm(\Gn)=\fr{1}{\pi i}\int_0^\infty\fr{f_j(\tau)}{H^+(\tau)}\fr{\tau d\tau}{\tau^2-\Gn^2}=O\left(\fr{1}{\Gn^2}\right), 
\quad \Gn\in {\Bbb C}^\pm, \quad \Gn\to\infty.
\label{2.31}
\eeq
Clearly, the functions $\GF_0^\pm(\Gn)$ and  $\GF_1^\pm(\Gn)$ have a simple zero at the infinite point as it is required.
Due to the presence of the functions $f_j(\Gn)$ they have four unknown constants $c_m$ ($m=0,1,2,3$).

Consider now case (ii)  when one of the zeros of the polynomial $q(\Gn)$, $z_0=\Gn_0$, is real,
and the other two, $z_1=\Gn_1$ and $z_2=-\Gn_2$, lie in the half-planes ${\Bbb C}^+$ and ${\Bbb C}^-$, respectively.
Owing to this, we select the Wiener-Hopf factors in (\ref{2.26}) as
\beq
H^+(\Gn)=\fr{\Gn+\Gn_2}{(\Gn+\Gn_0)(\Gn+\Gn_1)}, 
\quad 
H^-(\Gn)=-\fr{\Gn-\Gn_2}{(\Gn-\Gn_0)(\Gn-\Gn_1)}, 
\label{2.62}
\eeq
Upon representing the function $f_j(\Gn)/H^+(\Gn)$ as the difference $\Psi^+_j(\Gn)-\Psi_j^-(\Gn)$
of the boundary values on the real axis of the Cauchy integral (\ref{2.28}) with the function $H^+(\tau)$ 
being the one in (\ref{2.62}) and using the asymptotics of $H^\pm(\Gn)=O(\Gn^{-1})$, $\Gn\to\infty$, we deduce
\beq
\GF^\pm_j(\Gn)=H^\pm(\Gn)[b_{j}+\Psi_j^\pm(\Gn)], \quad \Gn\in {\Bbb C}^\pm, \quad j=0,1.
\label{2.64}
\eeq
Here, $b_0$ and $b_1$ are arbitrary constants. Because of the simple poles of the factors $H^\pm(\Gn)$ at 
$\mp\Gn_0$ in the real axis, the functions $\GF^\pm_j(\Gn)$ are not continuous at these points.
Due to the symmetry, to remove these singularities, it is necessary and sufficient to put 
\beq
b_{j}=-\Psi_j^+(-\Gn_0), \quad j=0,1.
\label{2.65}
\eeq
Because $f_j(\Gn)/H^+(\Gn)=0$
at $\Gn=-\Gn_0$ the condition (\ref{2.65}) is equivalent to
\beq
b_{j}=-\fr{1}{2\pi i}\int_{-\infty}^\infty
\fr{f_j(\tau)(\tau+\Gn_1)d\tau}{\tau+\Gn_2}, \quad j=0,1.
\label{2.66}
\eeq
Thus, as in case (i), the solution of the problem has only four arbitrary constants of the functions
$f_j(\Gn)$, $c_0,\ldots,c_3$, and $\GF_j^\pm(\Gn)=O(\Gn^{-1})$, $\Gn\to\infty$.

In case (iii), two zeros  of the polynomial $q(\Gn)$ lie in the upper half-plane,
$z_0=\Gn_0$ and $z_1=\Gn_1$, and the third one lies in the lower half-plane, $z_2=-\Gn_2$. The index of the 
function $H(\Gn)$ is now equal to 1.
We employ the same factorization as in the previous case. However, the factors $H^+(\Gn)$
and $H^-(\Gn)$ given by (\ref{2.62}) do not have poles on the real axis. What
is common with case (ii) is the asymptotics of the factors at the infinite point, $H^\pm(\Gn)=O(\Gn^{-1})$,
$\Gn\to\infty$. Therefore the solution of each Riemann-Hilbert problems (\ref{2.23}), (\ref{2.23'})
has an arbitrary constant and has the form 
\beq
\GF^\pm_j(\Gn)=H^\pm(\Gn)[b_{j}+\Psi_j^\pm(\Gn)],  \quad \Gn\in {\Bbb C}^\pm, \quad b_j=\const, \quad  j=0,1. 
\label{2.66.0}
\eeq
On the contrary to the previous case, there is no additional
condition (\ref{2.66}), and the constants $b_0$ and $b_1$ remain undetermined.

\subsection{Determination of the unknown constants $c_j$ ($j=1,2,3,4$)}\label{const}

For simplicity, we assume that the three functions $g_j$ ($j=0,1,2$) in the boundary conditions (\ref{2.2}) vanish,
$g_0(x)=g_1(x)=0$, $0<x<\infty$, and $g_2(y)=0$, $0<y<a$,
and select the function $g(x,y)$ in (\ref{2.1}) as $g(x,y)=-\Gd(x-x^\circ)\Gd(y-y^\circ)$, $0<x^\circ<\infty$ and
$0<y^\circ<a$. Then 
the functions $f_j(\Gn)$ in the Riemann-Hilbert problems (\ref{2.23})
can be represented as
\beq
f_j(\Gn)=-\fr{1}{q(-\Gn)\GD(\Gn)}\left[\sum_{m=0}^3 c_mf_{j}^m(\Gn)+f_{j}^4(\Gn)\right],\quad j=0,1,
\label{2.32}
\eeq
where
$$
f_{j}^{j}(\Gn)=2(-1)^{j+1}\Gn[\mu_{1-j}\sinh a\Gz+(\Ga_{1-j}^2-\Gn^2)\Gz\cosh a\Gz](\Ga_2^2+\Gz^2),
$$
$$
f_j^{1-j}(\Gn)=2(-1)^{j}\Gn\Gz(\Ga_{j}^2-\Gn^2)(\Ga_2^2+\Gz^2),
$$
$$
f_j^{j+2}(\Gn)=2(-1)^{j+1}\Gn(\Ga_{j}^2-\Gn^2)[\mu_{1-j}\sinh a\Gz+(\Ga_{1-j}^2-\Gn^2)\Gz\cosh a\Gz],
$$
$$
f_j^{3-j}(\Gn)=2(-1)^{j}\Gn\Gz(\Ga_0^2-\Gn^2)(\Ga_1^2-\Gn^2), 
$$
$$
f_j^4(\Gn)=
[\Gn(\Ga_2^2+\Gz^2)\cos\Gn x^\circ+\mu_2\sin\Gn x^\circ]f_j^*(\Gn),\quad j=0,1,
$$
$$
f_0^*=-2(\Ga_0^2-\Gn^2)[\mu_1\sinh(y^\circ-a)\Gz-\Gz(\Ga_1^2-\Gn^2)\cosh(y^\circ-a)\Gz],
$$
\beq
f_1^*(\Gn)=2(\Ga_1^2-\Gn^2)[\mu_0\sinh y^\circ\Gz+\Gz(\Ga_0^2-\Gn^2)\cosh y^\circ\Gz],
\label{2.33}
\eeq
and 
\beq
\GD(\Gn)=[\Ga_1^2\mu_0+\Ga_0^2\mu_1-(\mu_0+\mu_1)\Gn^2]\Gz\cosh a\Gz+[\mu_0\mu_1+(\Ga_0^2-\Gn^2)(\Ga_1^2-\Gn^2)\Gz^2]\sinh a\Gz.
\label{2.34}
\eeq
The first two conditions for the constants $c_j$ ($j=0,\ldots,3$) come from the boundary conditions
(\ref{2.6}) which, in the case of consideration, may be written as
\beq
(-1)^{j+1}\tilde u_y(\Gn,y_{j})+\fr{\mu_j\tilde u(\Gn,y_{j})}{\Ga_j^2-\Gn^2}=\fr{(-1)^{j+1}c_j}{\Ga_j^2-\Gn^2},\quad j=0,1,
\label{2.35}
\eeq
where, as in (\ref{2.3}), $y_0=0$ and $y_1=a$. We invert the Laplace transforms in (\ref{2.35}) and
assert that  $\tilde u(\Gn,y_j)=\GF_{j}^+(\Gn)$, $j=0,1$. Due to the first two boundary condition in (\ref{2.5})
\beq
\lim_{x\to 0^+}\int_{-\infty}^\infty \tilde u_y(\Gn,y_j)e^{-i\Gn x}d\Gn=0,
\label{2.37}
\eeq
and one deduces from (\ref{2.35})
\beq
\int_{-\infty}^\infty \fr{[\mu_j\GF_{j}^+(\Gn)+(-1)^j c_j ]d\Gn}
{\Ga_j^2-\Gn^2}=0, 
\quad j=0,1,
\label{2.38}
\eeq
or, equivalently, 
\beq
\mu_j\GF_{j}^+(\Ga_j)+(-1)^{j}c_j=0, \quad j=0,1.
\label{2.39}
\eeq
 In view of formulas (\ref{2.29}), (\ref{2.64}), (\ref{2.66.0}) and the representation (\ref{2.32}) we can write down the following two equations
for the constants:
\beq
(-1)^jc_j+\mu_j H^+(\Ga_j)\left[\sum_{m=0}^3 c_m\psi_{j}^m(\Ga_j)+\psi^4_{j}(\Ga_j)+b_j\right]=0, \quad j=0,1,
\label{2.40}
\eeq
where $b_j=0$ in case (i),
\beq
b_j=-\sum_{m=0}^3 c_m\psi_{j}^m(-\Gn_0)-\psi^4_{j}(-\Gn_0)
\label{2.40.0}
\eeq
in case (ii), and $b_j$ are free constants in case (iii), and
\beq
\psi_{j}^m(\Gn)=-\fr{1}{2\pi i}\int_{-\infty}^\infty \fr{f_{j}^m(\tau)d\tau}{q(-\tau)\GD(\tau)H^+(\tau)(\tau-\Gn)},
\quad j=0,1,\quad m=0,\ldots,4.
\label{2.41}
\eeq
Analysis of the functions (\ref{2.33}) and (\ref{2.34}) shows that 
$$
|f_j^{m}(\Gn)|\le A_{jm}|\Gn|^{-3}, \quad \Gn\to\infty, \quad \Gn\in {\Bbb C},\quad A_{jm}=\const,\quad  j=0,1,\quad m=0,1,2,3,
$$
\beq
f_1^4(\Gn)=O(\Gn^{-1}e^{-y^\circ|\Gn|}),\quad 
f_2^4(\Gn)=O(\Gn^{-1}e^{-(a-y^\circ)|\Gn|}), \quad \Gn\to\pm\infty.
\label{2.41.0}
\eeq
Since $\I\Ga_j>0$, the Cauchy integrals (\ref{2.41}) in (\ref{2.40}) are not singular and can be evaluated
by simple numerical methods. Alternatively,  if $m\ne 4$, then series expansions of the integrals (\ref{2.41})  can by 
derived by the Cauchy residue theorem.
For  this approach, in case (i), we invoke the representation
\beq
\fr{f_j(\Gn)}{H^+(\Gn)}=\fr{1}{(\Gn^2-\Gn_0^2)(\Gn^2-\Gn_2^2)\GD(\Gn)}
\left[\sum_{m=0}^3 c_m f_{j}^m(\Gn)+ f_{j}^4(\Gn)\right],\quad j=0,1,
\label{2.41.1}
\eeq
that follows from (\ref{2.25.0}), (\ref{2.27}) and (\ref{2.32}).  Since $f_j^m(\Gn)/\GD(\Gn)$ is a meromorphic
function of $\Gn$, for $\Gn\in {\Bbb  C}^+$  this enables us to write 
$$
\psi_j^m(\Gn)=\fr{1}{2(\Gn_0^2-\Gn_2^2)}\left(-\fr{ f_{j}^m(-\Gn_0)}{\Gn_0\GD(-\Gn_0)(\Gn_0+\Gn)}
+\fr{ f_{j}^m(-\Gn_2)}{\Gn_2\GD(-\Gn_2)(\Gn_2+\Gn)}\right)
$$
\beq
+\sum_{s=0}^\infty\fr{f_{j}^m(-\tau_{s})}{(\tau_{s}^2-\Gn_0^2)(\tau_{s}^2-\Gn_2^2)\GD'(-\tau_{s})
(\tau_{s}+\Gn)},\quad m=0,1,2,3,
\label{2.41.3}
\eeq
and $\tau_{s}$ ($s=0,1,\ldots,$) are the zeros (they are all simple) of the functions $\GD(\Gn)/\Gz$ in the upper half-plane.
An analog of the series representation for $\Gn\in {\Bbb C}^+$  in cases (ii) and (iii) has the form
\beq
\psi_j^m(\Gn)=-\fr 
{f_j^m(-\Gn_2)}{2\Gn_2\GD(-\Gn_2)(\Gn_2+\Gn)}
+\sum_{s=0}^\infty\fr{f_{j}^m(-\tau_{s})}{(\tau_{s}^2-\Gn_2^2)\GD'(-\tau_{s})
(\tau_{s}+\Gn)}.
\label{2.41.4}
\eeq

We turn now to the inverse Laplace transform of the conditions  (\ref{2.19}) as $y\to 0^+$
and $y\to a^-$ and derive two more equations for the constants $c_j$ ($j=0,\ldots,3$).
Similarly to the previous case of the boundary conditions on the horizontal walls because of (\ref{2.5}) we have
\beq
u_x(0,y)=\fr{1}{2\pi i}\int_{-i\infty}^{i\infty} \hat u_x(0,i\Gz)e^{\Gz y}d\Gz\to 0, \quad y\to 0^+,\quad y\to a^-.
\label{2.42}
\eeq
This brings us two more equations for the constants $c_j$
\beq
-\mu_2\CJ(y_j)+\fr{c_3e^{i\Ga_2(a-y_j)}-c_2e^{i\Ga_2y_j}}{2i\Ga_2}=0, \quad j=0,1,
\label{2.44.1}
\eeq
where
\beq
\CJ(y)=\fr{1}{2\pi i}\int_{-i\infty}^{i\infty} \fr{\hat u(0,i\Gz)e^{\Gz y}d\Gz}{\Gz^2+\Ga_2^2}.
\label{2.44.2}
\eeq
According to (\ref{2.21}) and (\ref{2.22}) the function $\hat u(0,i\Gz)$ is given by
$$
\hat u(0,i\Gz)=-\fr{\sin\Gn x^\circ}{\Gn}e^{-y^\circ\Gz}-\fr{\GF_0^+(\Gn)-\GF_0^-(\Gn)}{2i\Gn(\Gn^2-\Ga_0^2)}[\mu_0-(\Gn^2-\Ga_0^2)\Gz]
$$
\beq
-\fr{\GF_1^+(\Gn)-\GF_1^-(\Gn)}{2i\Gn(\Gn^2-\Ga_1^2)}e^{-a\Gz}[\mu_1+(\Gn^2-\Ga_1^2)\Gz].
\label{2.46}
\eeq
The main difficulty in computing the integral in (\ref{2.44.2}) is the presence of the two-valued function 
$\Gn^2=\Gz^2+k^2$ in $\hat u(0,i\Gz)$. We make the substitution $\Gz=-i\Gx$ and fix a branch of the function
 $\Gn^2=k^2-\Gx^2$
in the $\Gx$-plane cut along the line joining the branch points $\Gx=\pm k$ and passing through the infinite point
$\Gx=\infty$.   Notice that due to the symmetry condition (\ref{2.23'}) the functions
\beq
\fr{\GF_j^+(\Gn)-\GF_j^-(\Gn)}{\Gn}=\fr{\GF_j^+(\sqrt{k^2-\Gx^2})-\GF_j^-(\sqrt{k^2-\Gx^2})}{\sqrt{k^2-\Gx^2}},
\quad j=0,1,
\label{2.47}
\eeq
are  meromorphic in the $\Gx$-plane and independent of the branch choice. On using the boundary condition of the Riemann-Hilbert problem
\beq
\GF_j^-(\Gn)=\fr{\GF_j^+(\Gn)-f_j(\Gn)}{H(\Gn)}, \quad j=0,1,
\label{2.48}
\eeq
one may continue analytically the functions $\GF_j^-(\Gn)$ into the upper $\Gn$-half-plane ${\Bbb C}^+_\Gn$
cut along the semi-infinite line $\{|\Gn|>k, \arg \Gn=\Ga\}$,  $\Ga=\arg k\in(0,\pi/2)$ and in a similar manner continue 
the functions $\GF_j^+(\Gn)$ from the half-plane ${\Bbb C}^+_\Gn$ into the lower $\Gn$-half-plane ${\Bbb C}^-_\Gn$ cut
along the ray  $\{|\Gn|>k, \arg \Gn=\Ga+\pi\}$. Because of the meromorphicity of
the functions (\ref{2.47}), the function $(\Gz^2+\Ga_2^2)^{-1}\hat u(0,i\Gz)$ is meromorphic everywhere
in the $\Gx$-plane,
and it
can equivalently be represented as
$$
\fr{\hat u(0,i\Gz)}{\Gz^2+\Ga_2^2}=\fr{i}{\Gn(\Ga_2^2+\Gz^2)+i\mu_2}\left\{
\fr{\mu_0-(\Gn^2-\Ga_0^2)\Gz}{\Gn^2-\Ga_0^2}\GF_0^+(\Gn)+
\fr{\mu_1+(\Gn^2-\Ga_1^2)\Gz}{\Gn^2-\Ga_1^2}e^{-a\Gz}\GF_1^+(\Gn)
\right.
$$
\beq
\left.
+e^{i\Gn x^\circ-y^\circ\Gz}-
\left(\fr{c_0}{\Ga_0^2-\Gn^2}+\fr{c_2}{\Gz^2+\Ga_2^2}\right)
+\left(\fr{c_1}{\Ga_1^2-\Gn^2}+\fr{c_3}{\Gz^2+\Ga_2^2}\right)e^{-a\Gz}\right\}, \quad \Gz=-i\Gx,
\label{2.49}
\eeq
for all $\Gx$ such that $\Gn\in {\Bbb C}^+_\Gn$. Similarly,
$$
\fr{\hat u(0,i\Gz)}{\Gz^2+\Ga_2^2}=\fr{i}{\Gn(\Ga_2^2+\Gz^2)-i\mu_2}\left\{
-\fr{\mu_0-(\Gn^2-\Ga_0^2)\Gz}{\Gn^2-\Ga_0^2}\GF_0^-(\Gn)-
\fr{\mu_1+(\Gn^2-\Ga_1^2)\Gz}{\Gn^2-\Ga_1^2}e^{-a\Gz}\GF_1^-(\Gn)
\right.
$$
\beq
\left.
-e^{-i\Gn x^\circ-y^\circ\Gz}+
\left(\fr{c_0}{\Ga_0^2-\Gn^2}+\fr{c_2}{\Gz^2+\Ga_2^2}\right)
-\left(\fr{c_1}{\Ga_1^2-\Gn^2}+\fr{c_3}{\Gz^2+\Ga_2^2}\right)e^{-a\Gz}\right\}, \quad \Gz=-i\Gx,
\label{2.50}
\eeq
when $\Gn$ lies in the lower $\Gn$-half-plane ${\Bbb C}^-_\Gn$.
Examine now the conditions on $\Gx$ which imply $\Gn\in {\Bbb C}^-_\Gn$ and $\Gn\in {\Bbb C}^+_\Gn$.
Fix the branch $\Gn=i\sqrt{\Gx^2-k^2}$  by the conditions
\beq
\Gx\mp k=\Gr_\pm e^{i\Gt_\pm}, \quad \Ga-2\pi<\Gt_+<\Ga, \quad \Ga-\pi<\Gt_-<\Ga+\pi.
\label{2.51}
\eeq
It will be convenient to split the $\Gx$-plane into the following six sectors:
$$
D_1^+=\{0<\arg\Gx<\Ga\}, \quad D_2^+=\{\Ga<\arg\Gx<\pi/2\},
$$$$
 D_3^+=\{\pi/2<\arg\Gx<\pi\},  \quad 
D_1^-=\{\pi<\arg\Gx<\Ga+\pi\}, 
$$
\beq
 D_2^+=\{\Ga+\pi<\arg\Gx<3\pi/2\},\quad D_3^+=\{3\pi/2<\arg\Gx<2\pi\}.
\label{2.52}
\eeq
On employing (\ref{2.51}) we discover that the branch $\Gn$ maps the 
sectors $D_1^\pm$ and $D_3^\pm$ into the upper $\Gn$-half-plane, while the sectors
$D_2^\pm$ are mapped into the lower half-plane, that is
\beq
\Gn: D_1^\pm\to {\Bbb C}^+_\Gn, \quad \Gn: D_2^\pm\to {\Bbb C}^-_\Gn, \quad 
\Gn: D_3^\pm\to {\Bbb C}^+_\Gn.
\label{2.53}
\eeq
Consider case (i). The functions $q(\Gn)=\Gn(\Ga_2^2+\Gz^2)+i\mu_2$  has three zeros in the $\Gn$-plane,
$-\Gn_0\in {\Bbb C}_\Gn^-$, $\Gn_1\in{\Bbb  C}_\Gn^+$, and  $-\Gn_2\in{\Bbb  C}_\Gn^-$.  Denote
$\Gx_1=i\sqrt{\Gn_1^2-k^2}\in {\Bbb C}^+$. The function $\hat u(0,i\Gz)/(\Gz^2+\Ga_2^2)$ has two 
poles $\pm\Gx_1\in {\Bbb C}^\pm$ associated with the zero $\Gn_1$.
Analysis of the integral  (\ref{2.44.2}) shows that regardless which formula (\ref{2.49})
or (\ref{2.50}) for the integrand  is used there is only one pole of
the functions $q(\Gn)$ and $q(-\Gn)$ that generates
a nonzero residue. In the former case this pole is $\Gn=\Gn_1\in {\Bbb C}_\Gn^+$, 
while in the case (\ref{2.50}) it
is $\Gn=-\Gn_1\in{\Bbb  C}_\Gn^-$, and because of the symmetry property (\ref{2.23'}) the results of the integration
are the same.

Recall that the zeros of the function $q(\Gn)$ in case (ii) are $\Gn_0$ ($\I\Gn_0=0$), $\Gn_1\in {\Bbb C}_\Gn^+$, and  $-\Gn_2\in {\Bbb C}_\Gn^-$, and in case (iii), $\Gn_0\in {\Bbb C}^+_\Gn$, $\Gn_1\in {\Bbb C}_\Gn^+$, and  $-\Gn_2\in {\Bbb C}_\Gn^-$.
Denote  $\Gx_l=i\sqrt{\Gn_l^2-k^2}\in {\Bbb C}^+$, $l=0,1$. Then the function $\hat u(0,i\Gz)/(\Gz^2+\Ga_2^2)$ has two 
poles $\pm\Gx_0\in {\Bbb C}^\pm$ associated with the zero $\Gn_0$ of $q(\Gn)$ and two 
poles $\pm\Gx_1\in {\Bbb C}^\pm$ due to the zero $\Gn_1$.

In all cases (i) to (iii), in addition to the poles  at the zeros of the functions  $\Gn(\Ga_2^2+\Gz^2)\pm i\mu_2$ in (\ref{2.49})
and (\ref{2.50}), respectively, the function $\hat u(0,i\Gz)/(\Gz^2+\Ga_2^2)$ possesses
simple poles $\hat\Gx_0=-i\sqrt{\Ga_0^2-k^2}\in {\Bbb C}^-$ and
$\hat\Gx_1=i\sqrt{\Ga_1^2-k^2}\in {\Bbb C}^+$. Select  $\Gn(\hat\Gx_j)=i\sqrt{\hat\Gx_j^2-\Ga_j^2}=\Ga_j$,
$j=0,1$. Finally, because of the presence of the function $1/(\Gz^2+\Ga_2^2)$ in both formulas, (\ref{2.49})
and (\ref{2.50}), the integrand has two poles $\Gx=\Ga_2$ and $\Gx=-\Ga_2$ lying on the upper and lower half-planes,
respectively.
By employing the Cauchy residue theorem 
we transform the integral (\ref{2.44.2}) to the form
$$
\CJ(y)=\sum_{m=s}^1\fr{1}{t_m}\left[-\left(\fr{c_0}{\Ga_0^2-\Gn_m^2}+
\fr{c_2}{\Ga_2^2-\Gx_m^2}\right)e^{i\Gx_m y}+\left(\fr{c_1}{\Ga_1^2-\Gn_m^2}+\fr{c_3}{\Ga_2^2-\Gx_m^2}\right)e^{i\Gx_m(a-y)}
\right.
$$
$$
\left.
+e^{i\Gx_m|y-y^\circ|+i\Gn_m x^\circ}+\Gr_{0m}\GF_0^+(\Gn_m)e^{i\Gx_my}+
\Gr_{1m}\GF_1^+(\Gn_m)e^{i\Gx_m(a-y)}\right]
$$
\beq
-\fr{c_0+\mu_0\GF_0^+(\Ga_0)}{2i\hat\Gx_0r_0}e^{-i\hat\Gx_0 y}-\fr{c_1-\mu_1\GF_1^+(\Ga_1)}{2i\hat\Gx_1r_1}e^{-i\hat\Gx_1 (a-y)}+\fr{-c_2e^{i\Ga_2 y}+c_3e^{i\Ga_2 (a-y)}}
{2i\mu_2 \Ga_2}.
\label{2.53.0}
\eeq
Here, $s=1$ in case (i) and $s=0$ in cases (ii), (iii),
$$
t_m=\fr{\Gx_m(\Ga_2^2-k^2+3\Gn_m^2)}{\Gn_m}, \quad 
 \Gr_{jm}=\fr{\mu_j+i\Gx_m(\Ga_j^2-\Gn_m^2)}{\Gn_m^2-\Ga_j^2}, 
$$
\beq
 r_j=-i\Ga_j(\Ga_2^2-\hat\Gx_j^2)+\mu_2,\quad j=0,1,\quad m=0,1.
\label{2.54}
\eeq
Owing to the relation (\ref{2.39}) we simplify the two equations (\ref{2.44.1}) to read
$$
\sum_{m=s}^1\fr{1}{t_m}\left[-\left(\fr{c_0}{\Ga_0^2-\Gn_m^2}+
\fr{c_2}{\Ga_2^2-\Gx_m^2}\right)e^{i\Gx_m y_j}+\left(\fr{c_1}{\Ga_1^2-\Gn_m^2}+\fr{c_3}{\Ga_2^2-\Gx_m^2}\right)e^{i\Gx_m(a-y_j)}
\right.
$$
\beq
\left.
+e^{i\Gx_m|y_j-y^\circ|+i\Gn_m x^\circ}+\Gr_{0m}\GF_0^+(\Gn_m)e^{i\Gx_m y_j}+
\Gr_{1m}\GF_1^+(\Gn_m)e^{i\Gx_m(a-y_j)}\right], \quad j=0,1.
\label{2.55}
\eeq
Now, upon plugging (\ref{2.29}), (\ref{2.64}) and (\ref{2.66.0})  in the relations (\ref{2.55}), we deduce the following
two
equations for the constants $c_j$ ($j=0,\ldots,3$):
\beq
D_{j0}c_0+D_{j1}c_1+D_{j2}c_2+D_{j3}c_3=E_j, \quad j=0,1,
\label{2.57}
\eeq
where
$$
D_{jn}=\sum_{m=s}^1[d_{jnm}^0+d_{jnm}^1], \quad n=0,1,2,3, 
$$$$
d_{j0m}^0=\fr{e^{i\Gx_m y_j}}{\Gn_m^2-\Ga_0^2}, \quad
d_{j1m}^0=-\fr{e^{i\Gx_m (a-y_j)}}{\Gn_m^2-\Ga_1^2}, 
\quad
d_{j2m}^0=\fr{e^{i\Gx_m y_j}}{\Gx_m^2-\Ga_2^2}, \quad 
d_{j3m}^0=-\fr{e^{i\Gx_m(a-y_j)}}{\Gx_m^2-\Ga_2^2}, 
$$
$$
d_{jnm}^1=H^+(\Gn_m)[\Gr_{0m}\psi_{1}^n(\Gn_m)e^{i\Gx_m y_j}+\Gr_{1m}\psi_{2}^n(\Gn_m)e^{i\Gx_m (a-y_j)}],
\quad m=s,1,
$$
$$
 E_j=-\sum_{m=s}^1\left\{
 e^{i\Gx_m |y_j-y^\circ|+i\Gn_m x^\circ}+H^+(\Gn_m)
\right.
$$
\beq
\left.
\times\left[\Gr_{0m}(\Gy_{1}^4(\Gn_m)+b_1)e^{i\Gx_my_j}+\Gr_{1m}(\Gy_{2}^4(\Gn_m)+b_2)
e^{i\Gx_m(a-y_j)}\right]\right\}, \quad j=0,1.
\label{2.58}
\eeq
Here, $s=1$ in case (i), $s=0$ in cases (ii), (iii), 
and $\psi_j^m(\Gn_1)$ are determined by the quadrature (\ref{2.41}) and by the series (\ref{2.41.3}), (\ref{2.41.4}).
Equations (\ref{2.40}) and (\ref{2.57}) comprise a system of four equations with respect to the four constants
$c_j$ ($j=0,\ldots, 3$). In case (iii), the constants $b_j$ are still not determined, while in the other two cases, they are fixed: $b_j=0$ in case (i)  and $b_j$ are given by (\ref{2.40.0}) in case (ii).

\subsection{Analysis of the solution}

To write down the function $u(x,y)$ for any internal point in the semi-strip, one needs to know the function or its normal derivative either on the vertical part of the boundary $W_2=\{x=0, 0<y<a\}$, or on both horizontal sides 
$W_j=\{0<x<\infty, y=y_j\}$, $j=0,1$. Then the solution can be constructed in a standard manner by the method of integral transforms.
In fact, on having solved the Riemann-Hilbert problems one can determine not only the Laplace transforms of the functions $u(x,y_j)$, but also the Laplace transforms $\tilde u_y(x,y_j)$ from the boundary conditions (\ref{2.6})
and the Laplace transforms of the functions $u(0,y)$ and $u_x(0,y)$ by employing the relations (\ref{2.21}). 
Upon inverting these Laplace transforms we will have the function and its normal derivative available on the whole
boundary of the semi-strip. Application of the 
Green formula for the Helmholtz operator yields an integral representation of the solution inside the domain.

We start with the function $u(0,y)$ for $0<y<a$. By inverting the Laplace transform we have
\beq
u(0,y)=\fr{1}{2\pi i}\int_{-i\infty}^{i\infty} \hat u(0,i\Gz)e^{\Gz y}d\Gz.
\label{2.67}
\eeq
On employing the theory of residues, similarly to the previous section, one deduces
 $$
u(0,y)=\sum_{m=s}^1\fr{\Ga_2^2-\Gx^2_m}{t_m}\left[\left(\Gr_{0m}\GF_0^+(\Gn_m)-\fr{c_0}{\Ga_0^2-\Gn_m^2}-
\fr{c_2}{\Ga_2^2-\Gx_m^2}\right)e^{i\Gx_m y}
\right.
$$
\beq
\left.
+\left(\Gr_{1m}\GF_1^+(\Gn_m)+\fr{c_1}{\Ga_1^2-\Gn_m^2}+\fr{c_3}{\Ga_2^2-\Gx_m^2}\right)e^{i\Gx_m(a-y)}
+e^{i\Gx_m|y-y^\circ|+i\Gn_m x^\circ}\right].
\label{2.68}
\eeq
Here, $s=1$ in case (i), $s=0$ in cases (ii), (iii), and 
\beq
\GF^+_j(\Gn_m)=H^+(\Gn_m)\left[\sum_{n=0}^3 c_n\psi_j^n(\Gn_m)+\psi_j^4(\Gn_m)+b_j\right], \quad j=0,1.
\label{2.69}
\eeq
Next we wish to determine the function $u$ on the two horizontal boundaries of the half-strip,
\beq
u(x,y_j)=\fr{1}{2\pi}\int_{-\infty}^\infty \GF_j^+(\Gn) e^{-\Gn x} d\Gn, \quad 0<x<\infty.
\label{2.70}
\eeq
On continuing analytically the functions $\GF_j^+(\Gn)$ into the lower half-plane by making use of the Riemann-Hilbert
boundary conditions (\ref{2.23}) we rewrite the representation (\ref{2.70}) as
\beq
u(x,y_j)=\fr{1}{2\pi}\int_{-\infty}^\infty \left\{H^+(\Gn)[\Psi^-(\Gn)+b_j]+f_j(\Gn)\right\} e^{-\Gn x} d\Gn, \quad 0<x<\infty.
\label{2.71}
\eeq
In general, the functions $u(0,y)$ and $u(x, y_j)$ derived do not satisfy the compatibility conditions
\beq
\lim_{x\to 0^+} u(x,0)=\lim_{y\to 0^+} u(0,y), \quad \lim_{x\to 0^+} u(x,a)=\lim_{y\to a^-} u(0,y), 
\label{2.72}
\eeq
which guarantee the continuity of the function $u(x,y)$ and therefore, due to  (\ref{2.3}), the continuity
of the pressure distribution $p(x,y)$ at the corners of the semi-strip.
In cases (i) and (ii), after the constants $c_j$  ($j=0,\ldots,3$) have been fixed by solving the system 
of four equations (\ref{2.40}) and  (\ref{2.57}), there is no way to satisfy the compatibility conditions
(\ref{2.72}), and  in general, both functions,  $u(x,y)$ and $p(x,y)$, are discontinuous at the corners.
The situation is different in case (iii). We still have two free constants $b_0$ and $b_1$. 
To meet the conditions (\ref{2.72}) in this case, we transform the integral (\ref{2.71}) by evaluating the weekly convergent part
$$
u(x,y_j)=-i\sum_{m=0}^1\fr{(-1)^m(\Gn_2-\Gn_m)}{\Gn_1-\Gn_0}[\Psi_j^-(-\Gn_m)+b_j]e^{i\Gn_m x}
$$
\beq
-\sum_{n=0}^3 c_n M_j^n(x)-M_j^4(x), \quad 0<x<\infty,
\label{2.73}
\eeq
where
$$
\Psi_j^-(-\Gn_m)=\sum_{n=0}^3 c_n\psi_j^n(-\Gn_m)+\psi_j^4(-\Gn_m),
$$
\beq
M_j^n(x)=\fr{1}{2\pi}\int_{-\infty}^\infty \fr{f_j^n(\Gn)e^{-i\Gn x}d\Gn}{q(-\Gn)\GD(\Gn)}, \quad n=0,\ldots,4.
\label{2.74}
\eeq
Furnished with the expressions
(\ref{2.68}) and (\ref{2.73}) of the function $u(x,y)$ on the boundary we are able to satisfy the conditions
(\ref{2.72}) and fix the remaining constants $b_0$ and $b_1$.  The two new equations have the form
\beq
\Gb_{j0} b_0+\Gb_{j1} b_1+\Gb b_j+\sum_{n=0}^3(\Gs_{jn}+\Gl_{jn})c_n=\nu_j, \quad j=0,1.
\label{2.75}
\eeq
Here, 
$$
\Gb_{j0}=-\sum_{m=0}^1\fr{(\Ga_2^2-\Gx_m^2)\Gr_{0m}}{t_m}H^+(\Gn_m)e^{i\Gx_m y_j},
\quad 
\Gb_{j1}=-\sum_{m=0}^1\fr{(\Ga_2^2-\Gx_m^2)\Gr_{1m}}{t_m}H^+(\Gn_m)e^{i\Gx_m (a-y_j)},
$$
$$
\Gb=-\fr{i}{\Gn_1-\Gn_0}\sum_{m=0}^1 (-1)^m(\Gn_2-\Gn_m),\quad
\Gl_{j0}=\sum_{m=0}^1\fr{(\Ga_2^2-\Gx_m^2)e^{i\Gx_m y_j}}{t_m(\Ga_0^2-\Gn_m^2)},
$$
$$
\Gl_{j1}=-\sum_{m=0}^1\fr{(\Ga_2^2-\Gx_m^2)e^{i\Gx_m (a-y_j)}}{t_m(\Ga_1^2-\Gn_m^2)},
\quad
\Gl_{j2}=\sum_{m=0}^1\fr{e^{i\Gx_m y_j}}{t_m},\quad
\Gl_{j3}=-\sum_{m=0}^1\fr{e^{i\Gx_m (a-y_j)}}{t_m},
$$
$$
\Gs_{jn}=-i\sum_{m=0}^1\fr{(-1)^m(\Gn_2-\Gn_m)}{\Gn_1-\Gn_0}[\psi_j^n(-\Gn_m)-M_j^n(0)
$$
$$
-\sum_{m=0}^1\fr{\Ga_2^2-\Gx_m^2}{t_m}H^+(\Gn_m)\left[\Gr_{0m}\psi_0^n(\Gn_m)e^{i\Gx_m y_j}+
\Gr_{1m}\psi_1^n(\Gn_m)e^{i\Gx_m (a-y_j)}\right],
$$
$$
\nu_j=M_j^4(0)+i\sum_{m=0}^1\fr{(-1)^m(\Gn_2-\Gn_m)}{\Gn_1-\Gn_0}\psi_j^4(-\Gn_m)+\sum_{m=0}^1
\fr{\Ga_2^2-\Gx_m^2}{t_m}
$$
\beq
\times\left\{
e^{i\Gx_m|y_j-y^\circ|-i\Gn_m x^\circ}+H^+(\Gn_m)\left[\Gr_{0m}\psi_0^4(\Gn_m)e^{i\Gx_m y_j}
+\Gr_{1m}\psi_1^4(\Gn_m)e^{i\Gx_m (a-y_j)}\right]\right\}.
\label{2.76}
\eeq
These equations combined with (\ref{2.40}) and (\ref{2.57}) form a system of six equations for
the six constants $c_0,\ldots, c_3$, $b_0$, and $b_1$. Therefore we may conclude (in general the matrix of the system
is not singular) that in case (iii) the solution of the boundary value problem (\ref{2.1}), (\ref{2.2}), (\ref{2.5}) 
exists, it is unique and satisfies the compatibility conditions (\ref{2.72}).

The results obtained are collected in the theorem below. 

Theorem 2.1. {\it Let $g(x,y)\in L_1({\Bbb R}^2)$, $g_j(x)\in L_1(0,\infty)$, $j=0,1$, $g_2(y)\in L_1(0,a)$ and let
these functions satisfy the Dirichlet conditions that is be piecewise monotonic and  have a finite number of discontinuities.

Suppose   $k=\Go/c$, $\Ga_j=\Go/\Gve_j$,  
$\Gm_j=\Go^2/\Gd_j$, where $c$, $\Gve_j$ and $\Gd_j$ are positive constants and $\Go=\Go_1+i\Go_2$, $\Go_j>0$,
$j=0,1$. Denote $y_0=0$ and  $y_1=a$.

Consider the boundary value problem
$$
(\GD+k^2)u(x,y)=g(x,y), \quad 0<x<\infty, \quad 0<y<a,
$$$$
u_{xxy}+\Ga_j^2 u_y-(-1)^j\mu_j u=g_j(x), \quad 0<x<\infty, \quad y=y_j, \quad j=0,1,
$$
\beq
u_{xyy}+\Ga_2^2u_x-\mu_2 u=g_2(y), \quad x=0, \quad 0<y<a,
\label{2.77}
\eeq
whose solution satisfies the four conditions
\beq
u_y(0^+,y_j)=0, \quad j=0,1; \quad u_x(0,0^+)=u_x(0, a^-)=0.
\label{2.78}
\eeq

Let the three zeros of the polynomial
$q(\Gn)=\Gn(\Gn^2-k^2+\Ga_2^2)+i\mu_2$ be $z_0$, $z_1$ and $z_2$. Then two zeros say, $z_1$ and $z_2$, lie
in the opposite half-planes, $\I z_1>0$ and $\I z_2<0$. For the third zero, $z_0$, there are three possibilities:
(i) $\I z_0<0$, (ii) $ \I z_0=0$, and (iii) $\I z_0>0$.

In all cases (i) to (iii) the solution of the problem  (\ref{2.77}) exists, and  the Dirichlet data
on the two horizontal sides of the semi-strip, $u(x,0)$ and $u(x,a)$, are expressed by (\ref{2.70}) 
through the solution of the two symmetric scalar Riemann-Hilbert problems (\ref{2.23}), (\ref{2.23'}), $\GF_0^+(\Gn)$
and $\GF_1^+(\Gn)$. 
In the first two cases these solutions given by (\ref{2.29}) and (\ref{2.64}), (\ref{2.66}),
respectively, have four arbitrary constants $c_j$ ($j=0,\ldots,3$). 
In case (iii), the functions  $\GF_0^+(\Gn)$
and $\GF_1^+(\Gn)$ have the form (\ref{2.66.0}) and possess six arbitrary constants  $c_j$ ($j=0,\ldots,3$), $b_0$, and $b_1$.

In particular, if $g_0(x)=g_1(x)=0$ $(0<x<\infty)$, $g_2(y)=0$ $(0<y<a)$, and $g(x,y)=-\Gd(x-x^\circ)\Gd(y-y^\circ)$,
$x^\circ\in (0,\infty)$, $y^\circ\in (0,a)$, then the edge conditions (\ref{2.78}) are equivalent to the system of
four linear algebraic
equations (\ref{2.40}), (\ref{2.57}), where $b_j=0$ in case (i), $b_j$ are given by (\ref{2.66}) in case (ii), and remain
free in case (iii). In general, in cases (i) and (ii), the function $u(x,y)$ is discontinuous at the edges $x=y=0$ and $x=0, y=a$. In case (iii), however, on fixing the constants $b_1$ and $b_2$ by solving the two equations 
(\ref{2.75}), it is possible to satisfy the compatibility conditions (\ref{2.72}) and find the unique solution
 of the problem (\ref{2.77}), (\ref{2.78}) continuous up to the boundary including the corners of the semi-strip. }

\setcounter{equation}{0}

\section{Semi-infinite waveguide: walls are elastic plates}\label{plates}

Our previous analysis of the Helmholtz equation (\ref{2.1}) in a semi-infinite strip has been entirely limited to the case
of membrane walls modeled by the third order boundary conditions (\ref{2.2}). We next turn to  the two-dimensional model problem 
of a compressible fluid bounded by elastic walls.  This brings us boundary conditions with derivatives of order five. Assume that $B_j$ and $m_j$ are the bending stiffness
and mass per unit area of the plate $W_j$, respectively ($j=0,1,2$). 
As in Section \ref{const}, the function $g(x,y)$ is taken to be $g(x,y)=-\Gd(x-x^\circ)\Gd(y-y^\circ)$,
$(x^\circ,y^\circ)$ is an internal point of the semi-infinite strip, and the governing equation is 
\beq
\left(\fr{\Md^2}{\Md x^2}+\fr{\Md^2}{\Md y^2}+k^2\right)u(x,y)=-\Gd(x-x^\circ)\Gd(y-y^\circ), 
\quad 0<x<\infty, \quad 0<y<a,
\label{3.1}
\eeq
For the walls modeled by thin elastic
plates under flexural vibrations the boundary conditions (\ref{2.2}) are replaced by (Leppington, 1978)
$$
\left[\left(\fr{\Md^4}{\Md x^4}-\Ga_0^4\right)\fr{\Md}{\Md y}+\mu_0\right]u=0, \quad  (x,y)\in W_0=\{0<x<\infty, y=0\},
$$
$$
\left[-\left(\fr{\Md^4}{\Md x^4}-\Ga_1^4\right)\fr{\Md}{\Md y}+\mu_1\right]u=0, \quad  (x,y)\in W_1=\{0<x<\infty, y=a\},
$$
\beq
\left[\left(\fr{\Md^4}{\Md y^4}-\Ga_2^4\right)\fr{\Md}{\Md x}+\mu_2\right]u=0, \quad  (x,y)\in W_2=\{x=0, 0<y<a\}.
\label{3.2}
\eeq
Here, $\Go=\Go_1+i\Go_2$, $\Go_j>0$, $\R[e^{-i\Go t} u(x,y)]$ is the fluid velocity potential introduced in Section
\ref{form}, 
\beq
\Ga_j^4=\fr{m_j\Go^2}{B_j}, \quad \mu_j=\fr{\Gr\Go^2}{B_j}, \quad j=0,1,2.
\label{3.3}
\eeq
We need to choose constraints at the two edges $x=0, y=0$ and $x=0, y=a$. It is  designated that  the plates
are clamped at the edges (Fig.1), and therefore the deflections and the angles of deflection equal zero at the edges,
$$
\fr{\Md u}{\Md y}(0^+,y_j)=\fr{\Md^2 u}{\Md x\Md y}(0^+,y_j)=0, \quad j=0,1,\quad y_0=0, \quad y_1=a,
$$
\beq
\fr{\Md u}{\Md x}(0,0^+)=\fr{\Md^2 u}{\Md x\Md y}(0,0^+)=0, \quad 
\fr{\Md u}{\Md x}(0,a^-)=\fr{\Md^2 u}{\Md x\Md y}(0,a^-)=0.
\label{3.4}
\eeq
On following the procedure presented in Section 2 we apply the Laplace transform (\ref{2.5'})
to the boundary value problem (\ref{3.1}) to (\ref{3.3}), integrate by parts and deduce the one-dimensional
boundary value problem (\ref{2.6}), where
$$
 f(y)=u_x(0,y)-i\Gn u(0,y)-e^{i\Gn x^\circ}\Gd(y-y_0),\quad \tilde\mu_j(\Gn)=\fr{\mu_j}{\Ga_j^4-\Gn^4},
$$
\beq
\tilde g^j(\Gn)=\fr{c_{j0}-i\Gn c_{j1}}{\Ga_j^4-\Gn^4}, \quad j=0,1,
\label{3.5}
\eeq
and
\beq
c_{j0}=\fr{\Md^4 u}{\Md x^3\Md y}(0^+,y_j), \quad c_{j1}=\fr{\Md^3 u}{\Md x^2\Md y}(0^+,y_j), \quad j=0,1.
\label{3.6}
\eeq
Since the only difference between the problem (\ref{2.6}) obtained in the previous section and the one 
derived here is the form of the functions $\tilde \mu_j(\Gn)$, $f(y)$ and $\tilde g^j(\Gn)$, we still
have the relations (\ref{2.16}) to (\ref{2.18}). The next step of the procedure of Section 2 is to apply
the Laplace transform to the boundary condition on the vertical wall, the third condition in (\ref{3.2}).
This brings us to equation (\ref{2.19}) with the following notations adopted for the problem under consideration:
\beq 
\hat\mu_2(\Gz)=\fr{\mu_2}{\Ga_2^4-\Gz^4},\quad 
\hat g^2(i\Gz)=\fr{c_{20}+\Gz c_{21}-(c_{30}+\Gz c_{31})e^{-a\Gz}}
{\Ga_2^4-\Gz^4},
\label{3.7}
\eeq
where
\beq
c_{20}=u_{xyyy}(0,0^+), \quad c_{21}=u_{xyy}(0,0^+),
\quad 
c_{30}=u_{xyyy}(0,a^-), \quad c_{31}=u_{xyy}(0,a^-).
\label{3.7.0}
\eeq
Analogously to Section 2 the functions $\GF_1^\pm(\Gn)=\tilde u(\pm\Gn, 0)$ and
 $\GF_2^\pm(\Gn)=\tilde u(\pm\Gn, a)$ solve the symmetric 
 Riemann-Hilbert problem (\ref{2.23}), (\ref{2.23'}) with the coefficient:
 \beq
 H(\Gn)=-\fr{\Gn+i\hat\mu_2(\Gz)}{\Gn-i\hat\mu_2(\Gz)}
 =-\fr{\Gn[\Ga_2^4-(\Gn^2-k^2)^2]+i\mu_2}{\Gn[\Ga_2^4-(\Gn^2-k^2)^2]-i\mu_2}.
 \label{3.8}
 \eeq
 Remarkably, the coefficient $H(\Gn)$ and its Wiener-Hopf factors share the main features of 
those derived for the membrane walls of the wave guide.
Let $Q(\Gn)=\Gn[(\Gn^2-k^2)^2-\Ga_2^4]-i\mu_2$, and $z_j$ ($j=0,1,\ldots,4$) be the zeros of this
polynomial. Then $H(\Gn)=Q(\Gn)/Q(-\Gn)$, and $\Gn=-z_j$ are the zeros of the
denominators in (\ref{3.8}). It turns out that for all realistic values of the problem parameters two zeros, $z_1=\Gn_1$ and $z_2=\Gn_2$,
lie in the upper half-plane ${\Bbb C}^+$, and two zeros, $z_3=-\Gn_3$ and $z_4=-\Gn_4$, are located
in the lower half-plane ${\Bbb C}^-$ ($\I\Gn_j>0$, $j=1,2,3,4$). As for the fifth zero, $z_0$, there are three possible cases,

(i) $z_0=-\Gn_0\in {\Bbb C}^-$, 

(ii) $z_0=\Gn_0\in{\Bbb R}$, and

(iii)  $z_0=\Gn_0\in {\Bbb C}^+$.

\vspace{.1in}

Table 2. The roots $z_j$ $(j=0,1,\ldots,4)$ of the polynomial $Q(\Gn)$ for the wave number $k=1+0.1i$ and some values of the parameters $\Gg_0$
and $\Gg_1$.

\vspace{1mm}
\begin{tabular}{|c|c|c|c|c|c|}
\hline
$z_j$ & $\Gg_0=5, \Gg_1=1$ & $\Gg_0=1, \Gg_1=0.1$ & $\Gg_0=1, \Gg_1=1$ & $\Gg_0=.1, \Gg_1=1$ \\
\hline
$z_0$& $-1.806 - 0.04917 i$  &$-1.414 - 0.09353 i$ & $-1.441 + 0.008144 i$ & $-1.319 + 0.1075i$\\
\hline
$z_1$&$ -0.02056 + 1.256 i$ & $0.08369 + 0.3690 i$              &$0.02245 + 0.7374 i$ & $0.02935 + 0.5892 i$\\
\hline 
$z_2$& $1.809 + 0.1846 i$ &    $1.416+ 0.1184 i$                 &  $1.448 + 0.2135 i $ & $1.320 + 0.2936 i $\\
\hline
$z_3$&$-0.09151 - 0.7698 i$ &    $-0.3550 - 0.2809 i$               &$-0.6625 - 0.5350 i $ & $-0.8586 - 0.5673i $\\
\hline
$z_4$&$ 0.1083 - 0.6219 i$     & $0.2701 - 0.1131 i$                     & $0.6330 - 0.4240 i $ & $0.8280 - 0.4229 i$\\
\hline
\end{tabular}

\vspace{.1in}

In Table 2, we show the roots of the polynomial $Q(\Gn)$ for some values of its parameters.
The following notations are adopted:
$\Ga_2^4=\Gg_0 k^2$ and $\Gm_2=\Gg_1 k^2$, where $\Gg_0=m_2c^2/B_2$ and $\Gg_1=\Gr c^2/B_2$.

In view of the properties of the zeros and poles of the function $H(\Gn)$ we split the function $H(\Gn)$ as $H(\Gn)=H^+(\Gn)/H^-(\Gn)$, $-\infty<\Gn<\infty$,
where
\beq
H^+(\Gn)=\fr{(\Gn+\Gn_0)(\Gn+\Gn_3)(\Gn+\Gn_4)}{(\Gn+\Gn_1)(\Gn+\Gn_2)}, \quad 
H^-(\Gn)=-\fr{(\Gn-\Gn_0)(\Gn-\Gn_3)(\Gn-\Gn_4)}{(\Gn-\Gn_1)(\Gn-\Gn_2)}
\label{3.9}
\eeq
in case (i) and
\beq
H^+(\Gn)=\fr{(\Gn+\Gn_3)(\Gn+\Gn_4)}{(\Gn+\Gn_0)(\Gn+\Gn_1)(\Gn+\Gn_2)}, \quad 
H^-(\Gn)=-\fr{(\Gn-\Gn_3)(\Gn-\Gn_4)}{(\Gn-\Gn_0)(\Gn-\Gn_1)(\Gn-\Gn_2)}
\label{3.10}
\eeq
in cases (ii) and (iii).

It is seen that the asymptotics of the factors $H^+(\Gn)$ and $H^-(\Gn)$ at infinity is the same as for the membrane
walls model and therefore the solution has the form 
\beq
\GF^\pm_j(\Gn)=H^\pm(\Gn)[b_{j}+\Psi_j^\pm(\Gn)],  \quad \Gn\in C^\pm, \quad  j=0,1. 
\label{3.11}
\eeq
where $\Psi(\Gn)$ is determined by (\ref{2.28}) and (\ref{2.31}),
$b_j=0$ in case (i), $b_j$ are expressed through $\Psi(-\Gn_0)$ by (\ref{2.65}) in case (ii) and 
$b_j$ are free constants in case (iii). Notice that now formula (\ref{2.65}) reads
\beq
b_{j}=-\fr{1}{2\pi i}\int_{-\infty}^\infty
\fr{f_j(\tau)(\tau+\Gn_1)(\tau+\Gn_2)d\tau}{(\tau+\Gn_3)(\tau+\Gn_4)}, \quad j=0,1.
\label{3.12}
\eeq 
The functions $f_j(\Gn)$ are given by (\ref{2.25}), (\ref{2.18}), where $\tilde\Gm_0$, $\tilde\Gm_1$, 
and $\hat\Gm_2$ have to be replaced by their expressions in (\ref{3.5}) and (\ref{3.7}).  The functions
$f_j(\Gn)$ possess eight free constants  $c_{j0}$ and $c_{j1}$, $j=0,1,2,3$, and for their determination
we have the same number of additional conditions (\ref{3.4}).  Similarly to Section \ref{const}
these edge conditions can be rewritten as a system of eight linear algebraic equations
for the eight constants $c_{j0}$ and $c_{j1}$.  The clamping edge conditions (\ref{3.4})
guarantee that the derivatives $u_y$ and $u_{yx}$ vanish when $x\to 0^+$ along the horizontal
walls $W_0$ and $W_1$, and the functions $u_y$ and $u_{xy}$ tend to zero as $y\to 0^+$ and $y\to a^-$
along the vertical wall $W_2$. As for the function $u(x,y)=(i\Go \Gr)^{-1} p(x,y)$, in general, it is discontinuous
in cases (i) and (ii). In the case (iii), as in Section \ref{const}, it is possible to achieve the continuity of the function 
$u(x,y)$ and therefore the continuity of the pressure distribution at the corners of the semi-strip.
This can be done by fixing the remaining free constants $b_1$ and $b_2$ on satisfying the compatibility
conditions (\ref{2.72}).  

\section{Conclusion}

We have developed  further the method of integral transforms and made it applicable to the Helmholtz 
equation in a semi-infinite strip $\{0<x<\infty, 0<y<a\}$ with higher order impedance boundary conditions.  
It has been shown that if the orders of the tangential 
derivatives in the functionals of the boundary conditions are  even numbers, then the problem
reduces to two symmetric scalar Riemann-Hilbert problems which share the same coefficient, $H(\Gn)$, and possess
different right-hand sides. 
The coefficient $H(\Gn)$ is  a rational function $P_n(\Gn)/P_n(-\Gn)$, where $n={\rm deg} P_n(\Gn)$, and
$n-1$ is the order of the tangential derivative on the side $\{x=0, 0<y<a\}$ of the semi-infinite strip. In the case $n=3$,
the corresponding boundary value problem for the Helmholtz equation models acoustic wave 
propagation in a semi-infinite waveguide whose walls are membranes, and if $n=5$, then the 
walls are elastic plates. It turns out that the right-hand sides of the Riemann-Hilbert problems associated with the membranes and elastic plates  possess four and eight free constants, respectively. 
We have shown how these constants can be fixed by the conditions at the two edges of the structure. 
It has been discovered that, in addition to these expected free constants, the solution may or may not have two more free constants. This depends on the index of the Riemann-Hilbert problems that in turn
is determined by the location of the zeros of the polynomial $P_n(\Gn)$ and, ultimately, 
by the three parameters, $k$, $\Gg_0$, and $\Gg_1$, where 
$k$ is the wave number, $\Gg_0=m_2c^2/T_2$, $\Gg_1=\Gr c^2/T_2$, $c$ is the sound speed in the fluid,
$\Gr$ is the mean fluid density, $m_2$ and $T_2$
are the mass per unit area and the surface tension (the membrane case), respectively, of the finite vertical wall of the semi-strip.
In the plate case $T_2$ is replaced by $B_2$, the bending stiffness of the plate $x=0, 0<y<a$.   
For acceptable values of the parameters, the index $\Gk$ of the symmetric Riemann-Hilbert problems is either
$-1$, and the solution is unique, or 1, and then each problem has its own free constant.
We have shown that if $\Gk=-1$, then the solution satisfies the edge conditions, but
the pressure distribution $p(x,y)$ is discontinuous at the two corners. In the case $\Gk=1$, the two free constants 
available can
be fixed such that the function $p(x,y)$ is continuous at the vertices $x=y=0$ and $x=0, y=a$.

The approach we have presented works if the  governing PDE is of order two, has only even order
derivatives, and the tangential derivatives in the generalized impedance boundary conditions
are of an even order. If the functional of the boundary conditions has  tangential derivatives 
of an odd order, as in the Poincar\'e boundary value problem, then the problem is transformed into an order-2  vector Riemann-Hilbert 
problem whose coefficient is explicitly factorized only in some particular cases.

\end{document}